\documentclass[a4paper]{article}
\usepackage[utf8]{inputenc}
\usepackage{color}
\usepackage{amsmath}
\usepackage{amsmath}
\usepackage{amssymb}
\usepackage{multirow}
\usepackage{multicol}
\usepackage{authblk}
\usepackage[normalem]{ulem}

\usepackage{graphicx}
\usepackage{color}
\usepackage[colorinlistoftodos]{todonotes}
\usepackage{enumerate}
\usepackage{tikz,xcolor}
\usetikzlibrary{patterns, calc, positioning, shapes.geometric, arrows, fit, arrows.meta, backgrounds}
\usepackage{float}
\usepackage{wrapfig}
\usepackage{subcaption}
\usepackage{caption}
\usepackage{mathtools}
\usepackage{pgfplots}
\pgfplotsset{compat=1.10}
\usepackage{xcolor}
\usepackage{amsthm}
\usepackage[margin=2.5cm]{geometry}
\usepackage{hyperref}
\newcommand{\R}{\mathbb{R}}
\newcommand{\dd}{\partial}
\newcommand{\dbar}{\overline{\partial}}
\newcommand{\DOm}{\partial\Omega}
\newcommand{\Om}{\Omega}

\newcommand{\C}{{\mathbb C}}

\newcommand{\bL}{\boldsymbol{\Lambda}}

\newcommand{\mbR}{\mathbf{R}}
\newcommand{\cR}{{\mathcal R}}
\newcommand*\dif{\mathop{}\!\mathrm{d}}

\newcommand{\vp}{\varphi}

\newcommand{\w}{v}
\newcommand{\tw}{\widehat \w}

\def \bfo {\begin {displaymath} }
\def \efo {\end {displaymath} }
\def \beq {\begin {eqnarray}}
\def \eeq {\end {eqnarray}}
\def \ba {\begin {eqnarray*}}
\def \ea {\end  {eqnarray*}}
\def \R {{\mathbb R}}

\def \Re{\text{Re}\,}
\def \re{\text{Re}\,}

\def\tilde{\widetilde}

\def\dbar{\overline\partial}

\def \mbeq {\begin {eqnarray}}
\def \meeq {\end {eqnarray}}

\def \bfo {\begin {displaymath} } 
\def \efo {\end {displaymath} } 

\def \beq {\begin {eqnarray}}
\def \eeq {\end {eqnarray}}
\def \ba {\begin {eqnarray*}}
\def \ea {\end  {eqnarray*}}

\def \F {{\cal F}}

\def \cR {{\cal R}} 
 
\def \C {{\mathbb {C}}} 
\def \R {{\mathbb {R}}}

\def \H2s {H^{s+1}_0(\partial \M\times [0,T/2])}

\def \e {\varepsilon}

\def \pa0 {\partial _0}

\def \re {\Re}

\def \p {\partial}

\def \tilde{\widetilde}

\def \e {\epsilon}
\def \A {\mathcal A}

\newcommand{\M}{{\mathcal M}}

\newcommand{\beqla}[1] {\begin {eqnarray}\label{#1}}

\def \e {\varepsilon}

\def \beq {\begin {eqnarray}}
\def \eeq {\end {eqnarray}}
\def \ba {\begin {eqnarray*}}
\def \ea {\end  {eqnarray*}}

\newtheorem{theorem}{Theorem}[section]

\tikzset{
    module/.style={%
        draw, rounded corners,
        minimum width=#1,
        minimum height=7mm,
        font=\sffamily
        },
    module/.default=2cm,
    >=LaTeX
}

\title{CT scans without X-rays: parallel-beam imaging from nonlinear current flows}

\author[a,1]{Melody Alsaker}
\author[b]{Siiri Rautio}
\author[c]{Fernando Moura}
\author[d]{Juan Pablo Agnelli}
\author[e]{Rashmi Murthy}
\author[b]{Matti Lassas}
\author[f]{Jennifer L.\ Mueller}
\author[b]{Samuli Siltanen}

\affil[a]{\small Department of Mathematics, Gonzaga University, Spokane, WA 99258, USA}
\affil[b]{Department of Mathematics and Statistics, University of Helsinki, Helsinki, Finland}
\affil[c]{Engineering, Modeling and Applied Social Sciences Center,
Federal University of ABC, São Paulo, Brazil}
\affil[d]{FaMAF, National University of Córdoba, Córdoba, Argentina and CIEM, National Scientific and Technical Research Council (CONICET)}
\affil[e]{Department of Mathematics, Bangalore University, Bangalore, India}
\affil[f]{Department of Mathematics \& School of Biomedical Engineering, Colorado State University, Fort Collins, CO 80521, USA}
\begin{document}
\date{}

\maketitle

\begin{abstract}
Parallel-beam X-ray computed tomography (CT) and electrical impedance tomography (EIT) are two imaging modalities which stem from completely different underlying physics, and for decades have been thought to have little in common either practically or mathematically. CT is only mildly ill-posed and uses straight X-rays as measurement energy,  which admits simple linear mathematics. However, CT relies on exposing targets to ionizing radiation and requires cumbersome  setups with expensive equipment. In contrast, EIT uses harmless electrical currents as measurement energy and can be implemented using simple low-cost portable setups. But EIT  is burdened by nonlinearity stemming from the curved paths of electrical currents, as well as extreme ill-posedness which causes characteristic low spatial resolution. In practical EIT reconstruction methods, nonlinearity and ill-posedness have been considered intertwined in a complicated fashion. In this work we demonstrate a surprising connection between CT and EIT which partly unravels the main problems of EIT and leads directly to a proposed imaging modality which we call virtual hybrid parallel-beam tomography (VHPT). We show that hidden deep within EIT data is information which possesses the same linear geometry as parallel-beam CT data. This admits a fundamental restructuring of EIT, separating ill-posedness and nonlinearity into simple modular sub-problems, and yields ``virtual radiographs'' and CT-like images which reveal previously concealed information. Furthermore, as proof of concept we present VHPT images of real-world objects.
\end{abstract}

\section*{Introduction}
One of the marvels of scientific investigation is the occasional discovery of a deep and unexpected connection between seemingly disparate phenomena. In the 19th century, for example, Maxwell surprised the scientific world by connecting electromagnetism and wave theory \cite{maxwell1865}, and R\"ontgen discovered by chance, while working with cathode rays, how electromagnetic waves in the form of X-rays can be used to view an object's hidden internal structures \cite{rontgen1896}. Analogously, we describe here the discovery of an unexpected link between nonlinear electrical imaging and parallel-beam straight-line tomography. This work establishes the existence of physical and mathematical connections between two very different imaging modalities: X-ray computed tomography (CT) and electrical impedance tomography (EIT), and describes the advances these connections bring to the fields of tomographic imaging and inverse problems.  

In CT imaging, a technology dating back to the 1970s \cite{hounsfield1973, ambrose1973}, straight X-ray beams pass through an object and the measured attenuation values are subjected to mathematical inversion, thus forming an image of internal structures.  We restrict our discussion of CT to 2D ``slice'' imaging using parallel-beam geometry, as used by first-generation CT scanners, in which we fix the direction of parallel X-ray beams and record the attenuation of each ray as it travels through the object. After logarithmic transformation a 1D profile function is produced, 
revealing the total attenuation along each ray. This measurement is repeated for a dense collection of directions, and profile functions from all angles are collected as the columns of a 2D projection array known as a sinogram. Since X-rays travel along straight lines, it is possible to derive a fairly simple linear geometric reconstruction method, as in the widely-used filtered back-projection (FBP) algorithm \cite{natterer1986computerized} based upon mathematics formulated by Radon in 1917 \cite{radon1917uber}, to transform the sinogram back to the spatial imaging domain. In the 21st century, iterative methods based on compressed sensing, such as total variation (TV) reconstruction \cite{rudin1992nonlinear}, have emerged to provide edge-preserving regularization in situations with limited raw data \cite{sidky2008image}. 

EIT on the other hand, which emerged in the 1980s \cite{brown1987}, is based on measuring the boundary effects of electrical current flows. 
We inject harmless electrical currents into a target and measure the  voltages arising on electrodes placed on the target's surface. This boundary data is used to produce an image of the internal electrical conductivity. The difficulties of EIT lie in its nonlinearity and extreme ill-posedness, which until the writing of this present work were considered to be intertwined in a complicated fashion. While X-rays travel in straight lines, this is not so for electrical current, which flows mostly through domains of highest conductivity. Thus in EIT the very structures we are attempting to image will themselves alter the travel paths of the measuring energy, leading to nonlinearity. Linear-ray reconstruction methods therefore cannot be directly applied.  Further, while image formation in both CT and EIT is ill-posed in the sense that different targets may produce almost the same measurements, CT is only mildly ill-posed \cite{natterer1986computerized} while EIT is extremely ill-posed \cite{alessandrini1988stable}, resulting in  reconstructions with characteristically low spatial resolution.  Our work demonstrates that EIT data can be nonlinearly processed into a sinogram having exactly the same parallel-beam geometry as data produced by first-generation CT scanners. This is remarkable, as curved electrical current flows differ drastically from the straight paths of X-rays.  

Among many EIT reconstruction approaches \cite{mueller2012linear,Martins2019}, we focus here on a class of methods based on exponentially-behaving complex geometric optics (CGO) solutions, which permit insight leading to a fundamental restructuring of the EIT problem as well as VHPT imaging.  We also restrict ourselves to imaging (approximately) 2D objects, but in practice it is common to produce ``slice'' images of 3D objects  by placing electrodes around the boundary of a cross-sectional plane, thus permitting real-world imaging scenarios.   

Our findings lead to both conceptual and practical advances presented in this work. Conceptually, we present an enhanced understanding of the extremely ill-posed nonlinear EIT problem.  This involves a modular decomposition where ill-posedness is confined to two \textit{linear} steps, and all remaining nonlinear steps allow mathematical interpretation. On the practical side, we propose an imaging modality which we call {\it virtual hybrid parallel-beam tomography} (VHPT), in which electrical measurements are processed into a sinogram which is produced ``virtually,'' without actual X-rays. We may then view the columns of the sinogram directly as ``virtual radiographs,'' or apply standard CT reconstruction techniques to produce a conductivity image. 

\section*{The Method}

We now provide an overview of our method, which invokes deep mathematical connections between EIT and CT, leading directly to practical VHPT imaging.

\begin{figure}[t]
\centering
  \includegraphics[width=\textwidth]{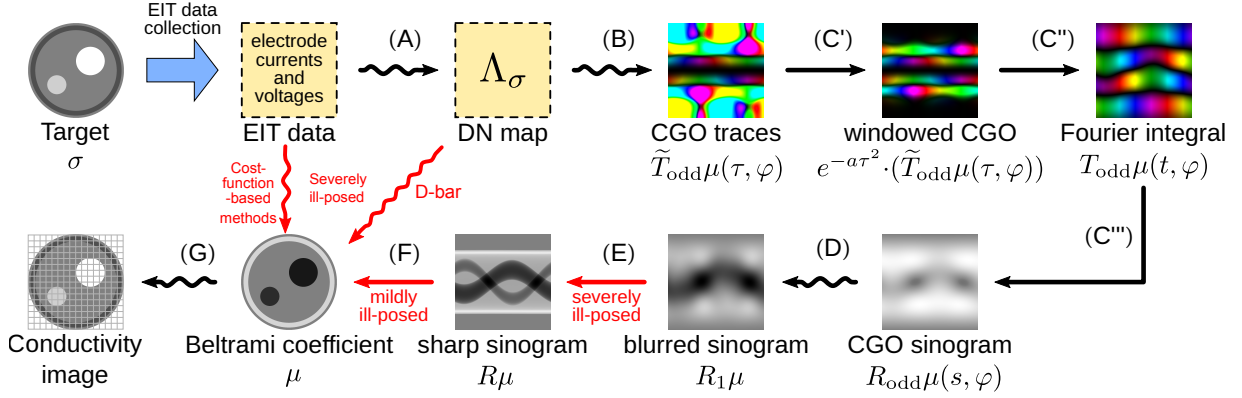}
  \caption{Method overview.  An illustration of our decomposition of the EIT problem into disjoint modules, with a representation of the virtual sinogram at bottom right. Arrow type and color indicate the type of operation; red: ill-posed, black: well-posed, straight: linear, and curved: nonlinear.  We flow through (A)-(G): (A) process EIT measurements into a DN map, (B) solve a boundary integral equation to yield CGO traces, (C$'$)-(C$'''$) apply a windowed Fourier transform and two simple linear operations, (D) remove higher-order scattering terms through machine learning, (E) apply 1D deconvolution separately to each column, (F) apply FBP, TV, or other CT reconstruction algorithm, (G) apply a simple algebraic formula. In this work, (D) is achieved through generalization of (E). Ill-posedness is completely confined to linear steps (E) and (F), with (E) containing almost all of it. Process flows for D-bar and cost-function-based methods are included for comparison.}
  \label{fig:VHEDprocess2}
\end{figure}
\subsection*{Mathematical Model}

Let $\Omega \subset \mathbb{R}^2$ be the unit disc and let the electrical conductivity $\sigma:\Omega \rightarrow \mathbb{R}$ satisfy
$ 0<c^{-1}\leq\sigma(x)\leq c, \; x \in \Omega.$ We begin by considering infinite-precision EIT measurements of the unknown conductivity $\sigma$ as described in \cite{Calderon1980}. In this ideal model called the {\it inverse conductivity problem}, one applies a voltage distribution along the boundary $\DOm$ and measures the current density through the boundary. Letting $u$ denote the electric potential,  the propagation of the electric field through $\Omega$, along with a known boundary voltage distribution, is modeled by the elliptic PDE
$\nabla\cdot\sigma\nabla u = 0$ in $\Om,$ with Dirichlet condition
$u|_{\DOm}$ specified. We seek to recover $\sigma$ satisfying this boundary value problem, given knowledge of the boundary voltage-to-current map, also called the Dirichlet-to-Neumann (DN) map $\Lambda_\sigma :  H^{1/2} (\DOm) \to H^{-1/2}(\DOm)$, given by $\Lambda_\sigma: u|_{\DOm} \longmapsto \sigma\frac{\partial u}{\partial \nu}  \; \mbox{on } {\DOm},$
where $\nu$ is the outward unit normal vector to $\partial \Omega$.

\subsection*{VHPT Imaging, Step-by-Step}

Here we describe the mathematical and computational steps of VHPT imaging,  following the diagram given in Figure \ref{fig:VHEDprocess2}. 

\subsubsection*{Step (A): EIT Data to DN Map} As a first step, we present an algorithm which transforms real-world data measurements into a mathematically ideal discrete DN map $\boldsymbol{\Lambda}_\sigma$.  This method may be viewed as a practical implementation of theory developed in  \cite{garde2021mimicking}. We assume two DN matrices computed from real-world EIT measurements: the ``target'' matrix $\boldsymbol{\Lambda}_{\mathrm{trg}}$ representing measurements collected at the boundary of the target object, and the ``calibration'' matrix $\boldsymbol{\Lambda}_{\mathrm{clb}}$ collected on a reference domain (e.g. a homogeneous tank). We further assume two simulated DN matrices corresponding to a homogeneous domain where $\sigma = 1$: $\boldsymbol{\Lambda}_1$, computed under the mathematically ideal continuum model  \cite{cheng1989};  and $\boldsymbol{\Lambda}_1^{\mathrm{cem}}$, simulated using the  ``complete electrode model'' \cite{cheng1989} which takes into account various physical effects of electrodes, by solving the forward EIT problem using finite element methods.  The final DN map is obtained via $\boldsymbol{\Lambda}_\sigma = (\tilde\bL_{\sigma}^{-1} + (\tilde\bL_{\sigma}^T)^{-1})/2$, where $\tilde\bL_{\sigma}= (\boldsymbol{\Lambda}_1^{\mathrm{cem}})^{-1} \boldsymbol{\Lambda}_{\mathrm{clb}} (\boldsymbol{\Lambda}_{\mathrm{trg}})^{-1} - (\boldsymbol{\Lambda}_1^{\mathrm{cem}})^{-1} + (\boldsymbol{\Lambda}_1)^{-1}.$ This is a form of difference imaging, as opposed to the absolute imaging scenario in which $\boldsymbol{\Lambda}_\sigma$ would be computed using only one set of experimentally collected measurements, $\boldsymbol{\Lambda}_\mathrm{trg}$. Difference imaging in EIT is well-known to provide significant advantages over absolute imaging in terms of robustness to noise and modeling errors \cite{brazey2022robust}. Details of this calibration approach are provided in the SI Appendix.

\subsubsection*{Step (B): DN Map to CGO Traces} Next, we follow  \cite{astala2006boundary,astala2011direct} to solve for the traces of certain complex geometric optics (CGO) functions $\omega^\pm = \omega^\pm(x,k)$ related to the conductivity $\sigma=\sigma(x)$, where $x$ is the spatial domain variable and $k$ is an introduced complex spectral parameter. 
Consider solutions $f_{\pm \mu}(x,k) = e^{ikx}(1 + \omega^{\pm}(x,k))$ of the Beltrami-type equation $\bar{\partial} f_\mu = \mu\overline{\partial f_\mu},$ where the Beltrami coefficient $\mu$ is given by $\mu :=\frac{1-\sigma}{1+\sigma},$
and $\partial =\frac 12(\frac{\partial}{\partial x_1}-i \frac{\partial}{\partial x_2})$ and 
$\overline{\partial }=\frac 12(\frac \partial {\partial x_1}+i\frac \partial {\partial x_2})$
are the complex Wirtinger derivatives. In~\cite{astala2006boundary} it is shown that $\omega^\pm|_{\DOm}$  may be obtained from knowledge of $\Lambda_\sigma$ by solving the boundary integral equation
\begin{equation} \label{eqn:BIE1}
{M_{\pm \mu}(\,\cdot\,,k)|_{\partial\Omega} + 1 = ({\mathcal P}_{\pm\mu}^k + {\mathcal P}_0)M_{\pm \mu}(\,\cdot\,,k)|_{\partial\Omega},}
\end{equation}
where $M_{\pm \mu}(x,k) := 1+\omega^{\pm}(x,k)$ and ${\mathcal P}_{\pm \mu}^k$ and ${\mathcal P}_0$ are certain projection operators computed directly from the DN map $\Lambda_\sigma$ (see the SI Appendix). It is further shown in~\cite{greenleaf2018propagation}  that the CGO solutions $\omega^{\pm}$ can be written as a scattering series
$\omega^{\pm} = \sum_{j=1}^{\infty} \omega_j^{\pm}.$
While the leading order term $\omega_1^+$ contains relevant information of the singularities of $\mu$ and thus of the conductivity $\sigma$, in the inverse problem setting with unknown conductivity we can only compute the full series $\omega^{\pm}$ and there is no formula to solve for $\omega_1^+$. However, in \cite{greenleaf2018propagation} it is also shown that by subtracting $\omega^-$ from $\omega^+$, one eliminates all even terms $\omega_{2j}^{+}$ in the scattering series of $\omega^+$. We therefore denote $\omega_{\text{odd}}^+ := \omega^+ - \omega^-$. Since in what follows we will focus on the ``plus'' cases $\omega_1^+$ and $\omega_{\text{odd}}^+$, for notational simplicity we denote these functions by $\omega_1$ and $\omega_{\text{odd}}$, respectively. After solving \eqref{eqn:BIE1} twice to obtain $\omega_{\text{odd}}$ from the solutions $M_{\pm \mu}$, we integrate the functions $\omega_{\text{odd}}$ to get a scattering operator $\widetilde{T}_{\text{odd}} \mu (\tau,\varphi)$ with parameters
$\tau \in \R$ and $\varphi \in [0,2\pi]$ (depicted in Figure~\ref{fig:VHEDprocess2}): 
\begin{equation*}
 \widetilde{T}_{\text{odd}} \mu (\tau,\varphi) := \frac{1}{2 \pi i} \int_{\partial\Omega} \omega_{\text{odd}}(x,\tau e^{i\varphi}) \dif x.
 \end{equation*}

\subsubsection*{Steps (C$'$)–(C$'''$): CGO Traces to CGO Sinogram} Noise in the data $\bL_\sigma$ causes inaccuracy in the solution of the boundary integral equation (\ref{eqn:BIE1}) which grows with $|\tau|$. As a low-pass filter, we apply a one-dimensional windowed Fourier transform to the scattering operator $\widetilde{T}_{\text{odd}} \mu (\tau,\varphi)$. We therefore restrict integration to an interval $[-R,R]$ and insert a Gaussian function $e^{-a\tau^2}$ with suitable parameter $a>0$ into the Fourier integral:
\begin{equation}\label{eq:T_odd}
T_{\text{odd}} \mu(t,\varphi) 
= \frac{1}{2 \pi i}\int_{-R}^{R} e^{-a\tau^2}e^{-i t\tau } (\widetilde{T}_{\text{odd}} \mu(\tau,\varphi))  \dif \tau,
\end{equation}
where the variable $t\in\mathbb{R}$ is called the {\it pseudo-time.} By this Fourier windowing, the noisy parts of the signal will be set to zero, but as a trade-off, the resulting  complex-valued $T_{\text{odd}} \mu(t,\varphi)$ has blurring in the pseudo-time domain. In practical computations, this corresponds to blurring along the columns of the matrix approximation to $T_{\text{odd}} \mu(t,\varphi)$, as depicted in Figure~\ref{fig:VHEDprocess2}. To compute the real-valued ``projection data'' function from the complex-valued version, we multiply with a complex phase factor and then integrate to get the CGO sinogram $
  R_{\text{odd}}\mu(s, \varphi)  = \int_{-\infty}^s \frac{e^{-i\varphi}}{2\pi i}  T_{\text{odd}} \mu(t, \varphi) \dif t.
$

\subsubsection*{Steps (D)–(E): CGO Sinogram to Sharp Sinogram}
The CGO sinogram is blurred due to the necessary Fourier windowing process.  Moreover,  due to our unavoidable use of $\omega_{\text{odd}}$ in place of the ideal $\omega_1$, the quantity $R_{\text{odd}}\mu(s, \varphi)$ does not exactly match  the desired theoretical framework. Based on \cite{greenleaf2018propagation}, the ideal quantity $\omega_1$ would correspond to the Fourier integral $T_1\mu$ and projection data $R_1\mu$, as opposed to $R_{\text{odd}}\mu$.  This discrepancy may be rectified by using a convolutional neural network (CNN) to remove the higher order terms of the scattering series in the projection domain, which corresponds to step (D) in Figure~\ref{fig:VHEDprocess2}. One could follow this machine-learning process with any linear deconvolution algorithm (step (E)) to correct the blur. In our implementation, however, we found it was efficient to perform steps (D) and (E) simultaneously by training a single neural network to perform the deconvolution operation $\text{NN}_{\theta}(R_1\mu) \approx R\mu,$
where $\theta$ is the learned parameter vector. The training inputs are ideal CGO sinograms computed from numerically simulated functions $\omega_1$, given by
$
R_1\mu(s,\varphi)= \int_{-\infty}^s \frac{-e^{i \varphi}}{2\pi i} (g*T_1)\mu(t,\varphi) \dif t,
$
where $g=g(t)$ is the inverse Fourier transform of $e^{-a\tau^2}$ and $T_1\mu(t,\varphi)$ is computed from $\omega_1$ and is analogous to~\eqref{eq:T_odd}. The  ground truth $R\mu$ is the Radon transform of $\mu$, defined as
$R\mu(s,\varphi) = \int_{x\cdot\varphi = s} \mu(x) dx.$
 After training, the network generalizes and can be used instead with inputs given by the CGO sinograms $R_{\text{odd}}\mu$ computed from the EIT data.
The network outputs sharp VHPT sinograms $S \approx R\mu$, where the higher-order scattering terms have consequently been removed. Details regarding network architecture and training are provided in the SI Appendix. 

\subsubsection*{Steps (F)–(G): Sharp Sinogram to Reconstruction} Remarkably, the sharpened sinogram embodies exactly the same geometry as one expects from parallel-beam X-ray tomography, but is produced without ionizing radiation. Geometric information about the target's interior may be extracted directly by viewing columns of the sinogram, but if desired we can transform this projection data into the spatial image domain using standard (mildly ill-posed) CT reconstruction methods. This linear inversion process (F) along with the previous linear deconvolution problem (E) together contain all ill-posedness of the VHPT chain, admitting an enhanced understanding of EIT within the framework of nonlinear ill-posed inverse problems.  Application of FBP, TV, or another standard CT reconstruction method to the sinogram will yield the Beltrami coefficient $\mu$. For this work we applied TV reconstruction, as described in greater detail in the SI appendix.  The conductivity $\sigma$  satisfies $\sigma = (1-\mu)/(1+\mu)$, and we apply this formula in step (G) to yield a conductivity image, thus completing the VHPT process.

\section*{Results} 

\subsection*{Conceptual Discoveries} Aside from providing a practical imaging technique, the previously described VHPT method yields a conceptual discovery which brings enhanced clarity to the EIT problem. The article \cite{greenleaf2018propagation} revealed a linear FBP framework hidden inside the nonlinear inverse problem of EIT, a breakthrough which was based on deep mathematical analysis \cite{duistermaat1972fourier, Greenleaf1989}. We now extend these results by presenting a fundamental restructuring of the decades-old EIT problem, as illustrated in Figure~\ref{fig:VHEDprocess2}. Here  we see the remarkable fact that the very complicated EIT problem may be modularly decomposed into a series of well-posed nonlinear problems (black arrows in Figure~\ref{fig:VHEDprocess2}), followed by two ill-posed but linear inverse problems (red arrows). One of these ill-posed modules is the standard (mildly ill-posed) CT reconstruction problem (solvable via FBP or TV, for example), and the other is simply a compilation of classical linear one-dimensional Gaussian deconvolution problems. In addition, all but one of the well-posed nonlinear problems are mathematically explicit. The sole non-explicit step in the process (step (D) in Figure~\ref{fig:VHEDprocess2}) can be viewed as stripping an infinite non-physical scattering series of its higher-order terms. In this work, we show that (D) and (E) may in fact be completed simultaneously using machine learning, although one can approach them as separate problems. The severely ill-posed and nonlinear processes labeled ``cost-function-based methods'' and ``D-bar'' in Figure~\ref{fig:VHEDprocess2} are not part of the VHPT imaging chain, but provide comparison with other  common EIT reconstruction algorithms wherein ill-posedness and nonlinearity are intertwined.

\subsection*{Practical Reconstructions}

\begin{figure} [t!]
  \raisebox{2cm}{I }\hfill\begin{subfigure}[b]{0.23\textwidth}
\includegraphics[width=\textwidth]{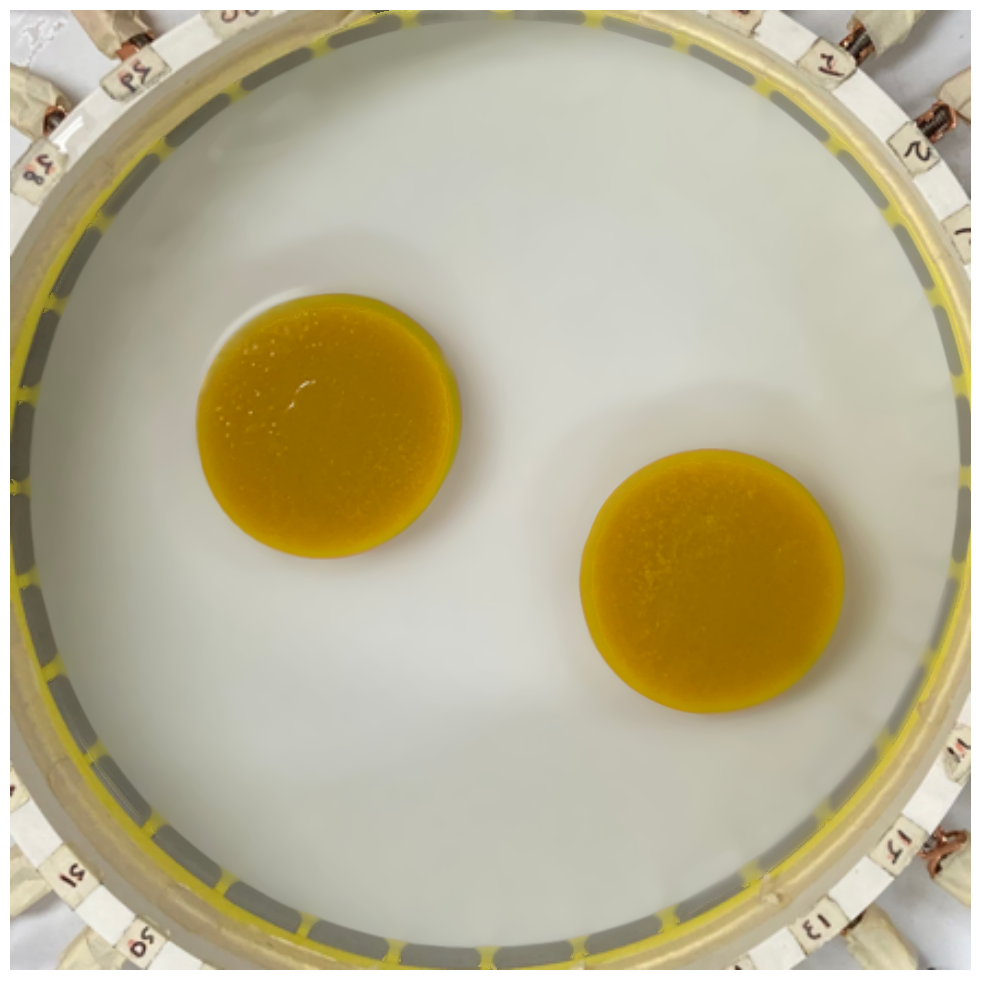}
    \label{fig:tank1}
  \end{subfigure}
  \hfill
  \begin{subfigure}[b]{0.44\textwidth}
    \includegraphics[width=\textwidth, height=3.7cm]{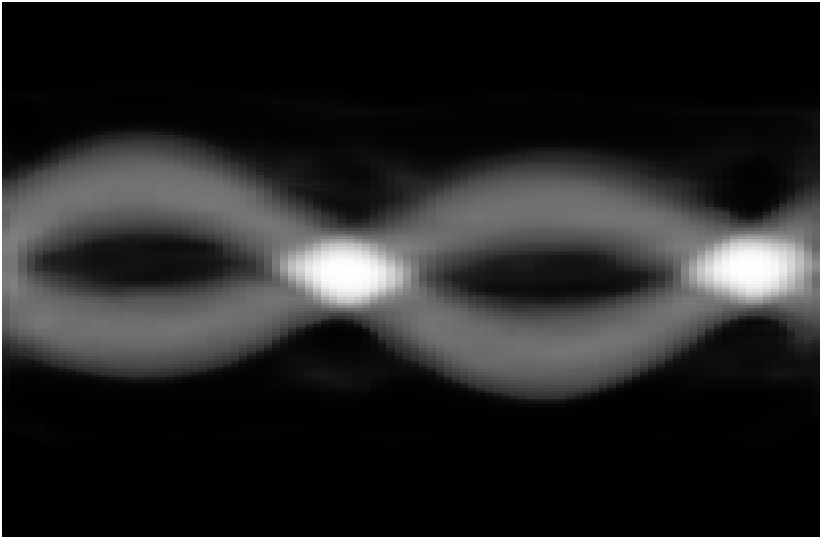}
    \label{fig:sinogram1}
  \end{subfigure}
  \hfill
  \begin{subfigure}[b]{0.23\textwidth}
    \includegraphics[width=\textwidth]{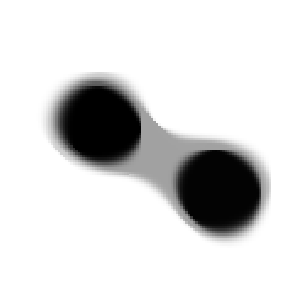}
    \label{fig:reco1}
  \end{subfigure}


    \raisebox{2cm}{II }\hfill\begin{subfigure}[b]{0.23\textwidth}
    \includegraphics[width=\textwidth]{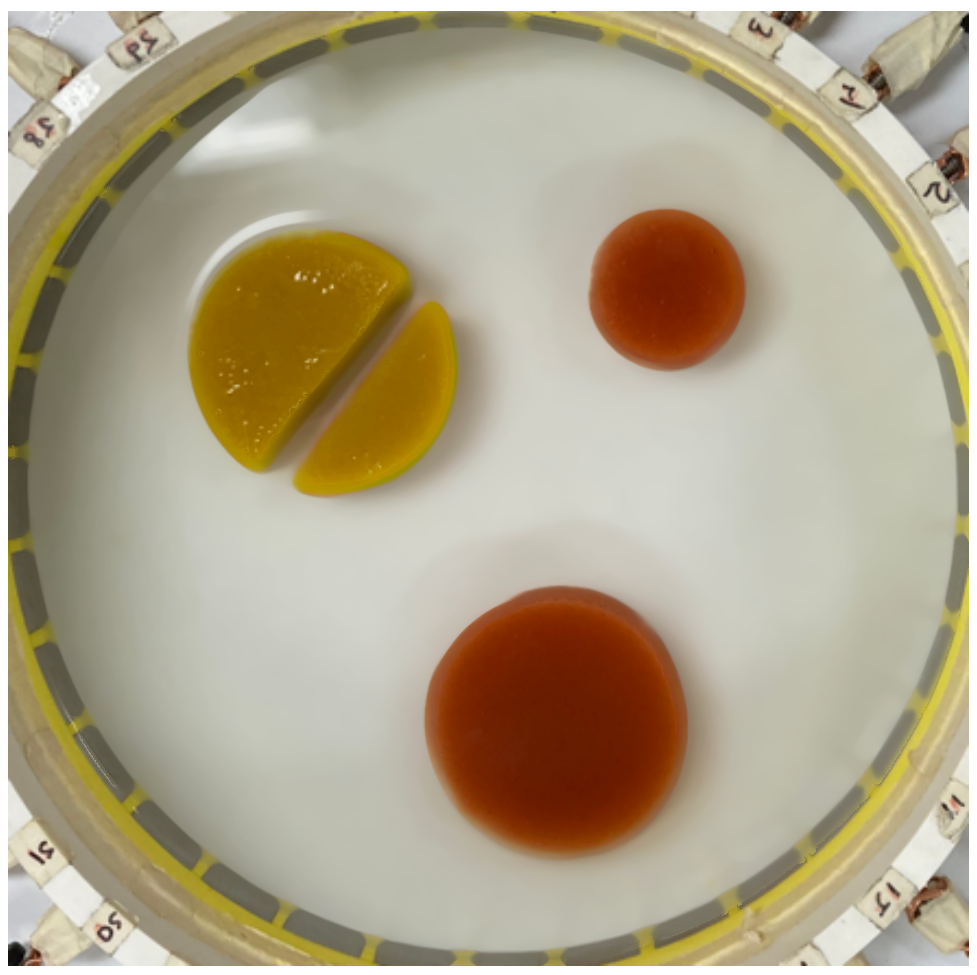}
    \label{fig:tank3}
  \end{subfigure}
  \hfill
  \begin{subfigure}[b]{0.44\textwidth}
    \includegraphics[width=\textwidth, height=3.7cm]{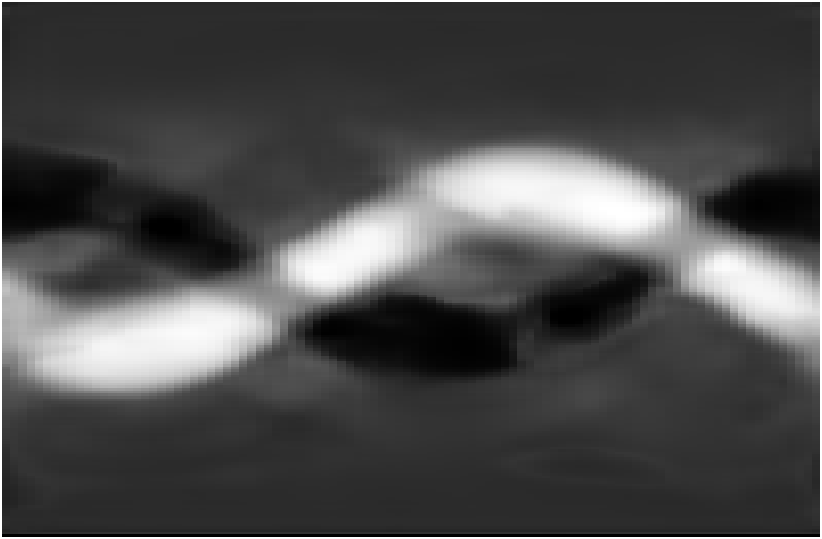}
    \label{fig:sinogram3}
  \end{subfigure}
  \hfill
  \begin{subfigure}[b]{0.23\textwidth}
    \includegraphics[width=\textwidth]{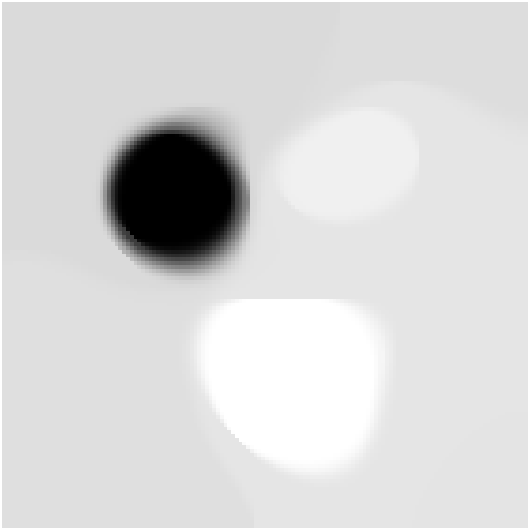}
    \label{fig:reco3}
  \end{subfigure}

    \raisebox{2cm}{III }\hfill \begin{subfigure}[b]{0.23\textwidth}
    \includegraphics[width=\textwidth]{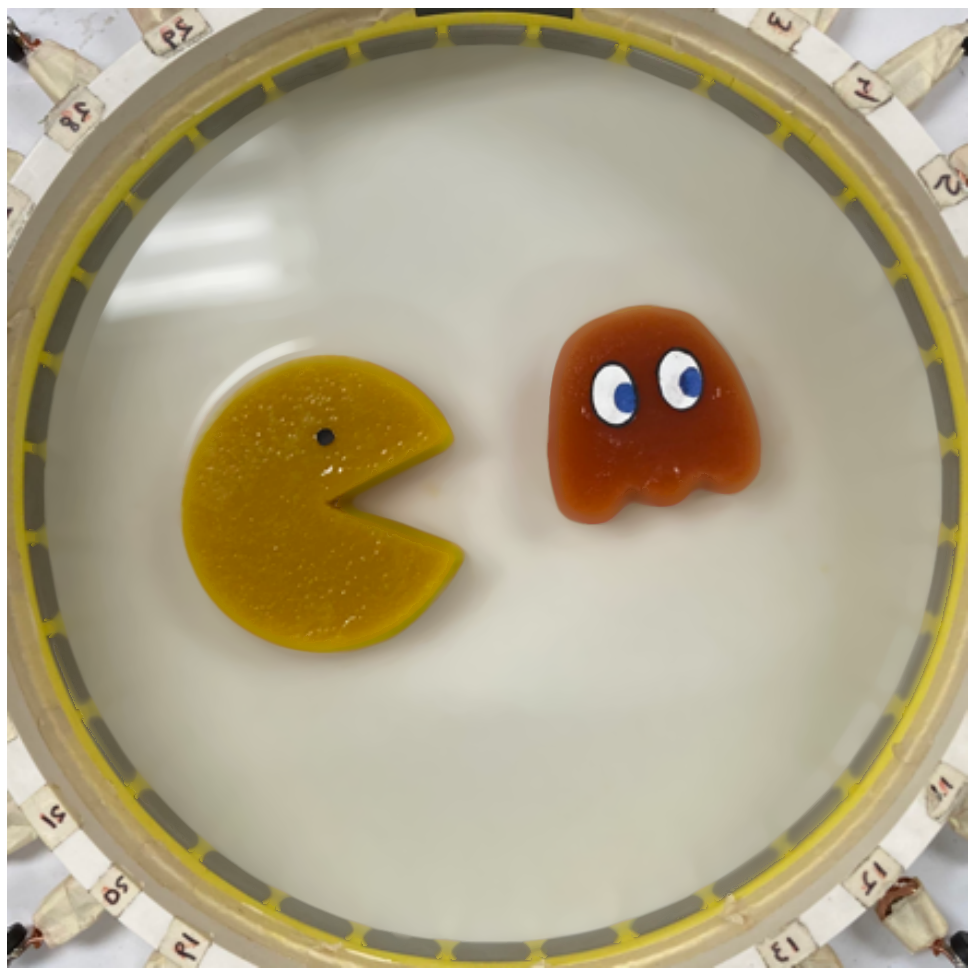}
    \caption*{A}
    \label{fig:tank4}
  \end{subfigure}
  \hfill
  \begin{subfigure}[b]{0.44\textwidth}
    \includegraphics[width=\textwidth, height=3.7cm]{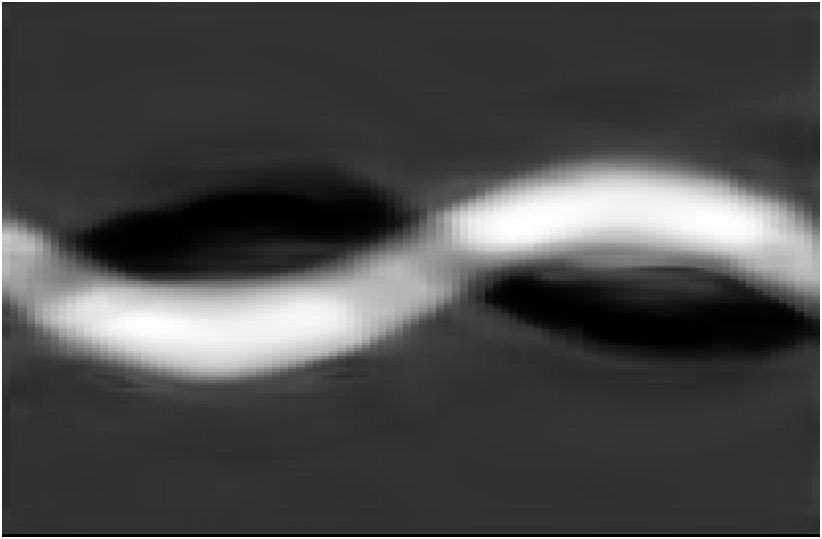}
    \caption*{B}
    \label{fig:sinogram4}
  \end{subfigure}
  \hfill
  \begin{subfigure}[b]{0.23\textwidth}
    \includegraphics[width=\textwidth]{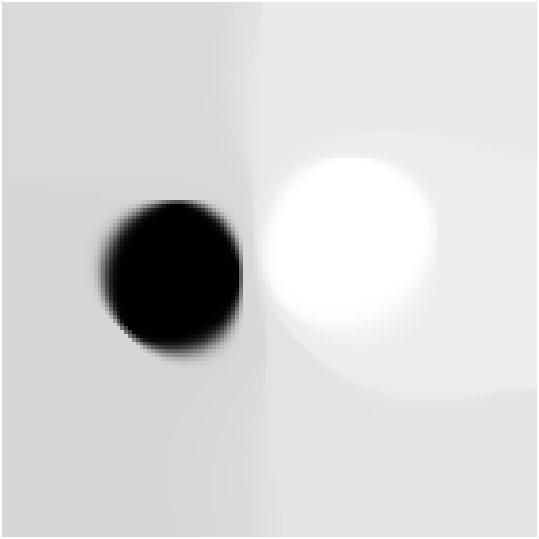}
    \caption*{C}
    \label{fig:reco4}
  \end{subfigure}
    \caption{Results. (A): The experimental setup for collection of the EIT target datasets I, II, III with 32 electrodes. Yellow agar inclusions are resistive and red inclusions are conductive, as compared to background. (B):  VHPT virtual (deblurred) sinogram for targets I, II, III. (C): VHPT conductivity reconstructions via total variation regularization for targets I, II, III. }
  \label{fig:results}
\end{figure}

In Figure~\ref{fig:results} we present VHPT conductivity reconstructions (C) of real-world objects, produced by applying total variation reconstruction \cite{bredies2014recovering} to virtual parallel-beam sinograms (B) computed from EIT measurements. Target datasets I-III, consisting of voltage measurements taken on all electrodes for each current injection, were collected using the ACT5 EIT system at Colorado State University \cite{ACT5paper}. Current injection at peak amplitude 0.35~mA was applied on 32 electrodes according to a trigonometric basis. We constructed each target to consist of a circular tank filled with saline solution of conductivity 2.7 mS/cm, with various homogeneous agar inclusions having conductivity either 0.5~mS/cm (yellow agar) or 3.7~mS/cm (red agar). 
The resulting VHPT reconstructions (C) in Figure~\ref{fig:results} depict conductivity contrasts of the target interior, ranging from lowest (black) to highest (white) conductivity as compared to background.  

\begin{figure}[t!]
  \raisebox{2.5cm}{I }\hfill\begin{subfigure}[b]{0.27\textwidth}
    \includegraphics[width=\textwidth]{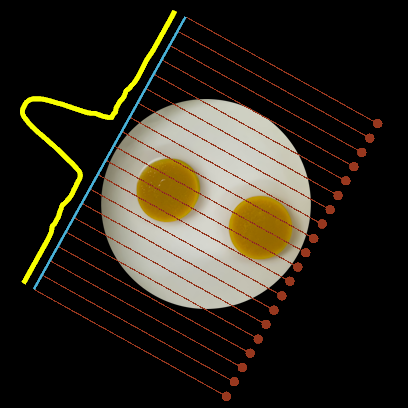}
    \label{fig:profile1a}
  \end{subfigure}
  \hfill
  \begin{subfigure}[b]{0.27\textwidth}
    \includegraphics[width=\textwidth]{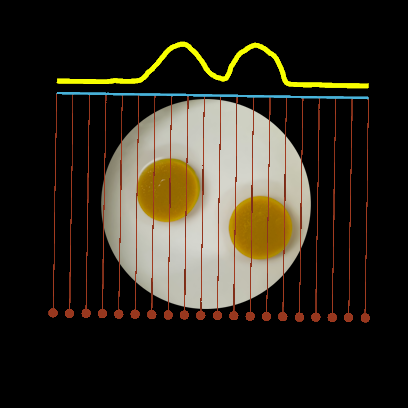}
    \label{fig:profile1b}
  \end{subfigure}
  \hfill
  \begin{subfigure}[b]{0.27\textwidth}
    \includegraphics[width=\textwidth]{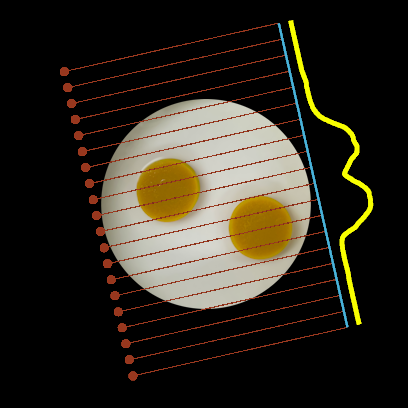}
    \label{fig:profile1c}
  \end{subfigure}

  
  \raisebox{2.5cm}{II }\hfill\begin{subfigure}[b]{0.27\textwidth}
    \includegraphics[width=\textwidth]{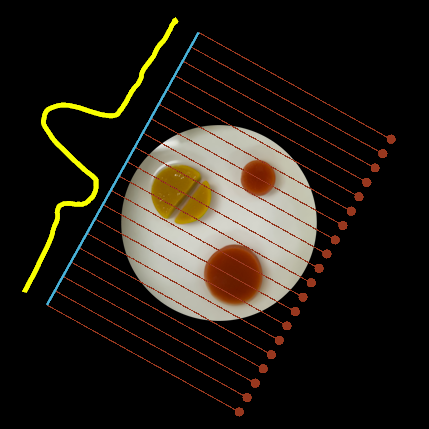}
    \label{fig:profile3a}
  \end{subfigure}
  \hfill
  \begin{subfigure}[b]{0.27\textwidth}
    \includegraphics[width=\textwidth]{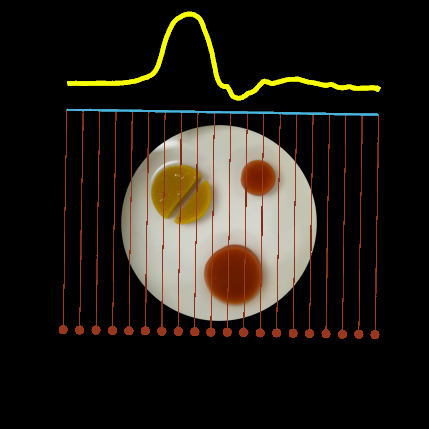}
    \label{fig:profile3b}
  \end{subfigure}
  \hfill
  \begin{subfigure}[b]{0.27\textwidth}
    \includegraphics[width=\textwidth]{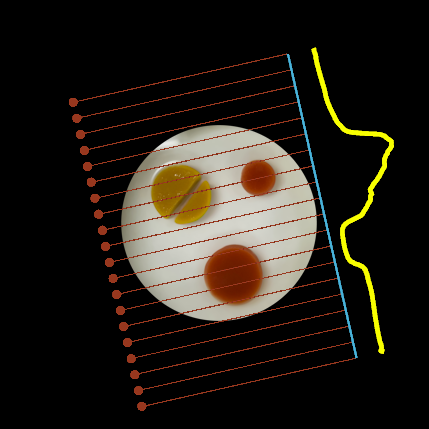}
    \label{fig:profile3c}
  \end{subfigure}

    \raisebox{2.5cm}{III }\hfill\begin{subfigure}[b]{0.27\textwidth}
    \includegraphics[width=\textwidth]{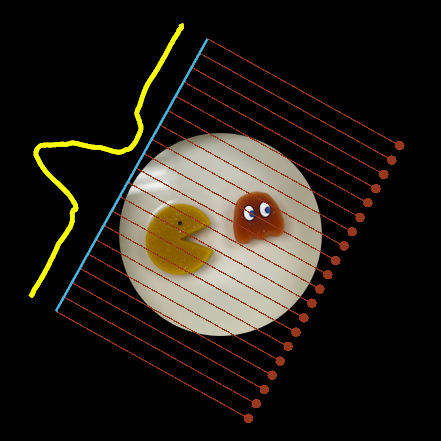}
    \caption*{A}
    \label{fig:profile4a}
  \end{subfigure}
    \hfill
  \begin{subfigure}[b]{0.27\textwidth}
    \includegraphics[width=\textwidth]{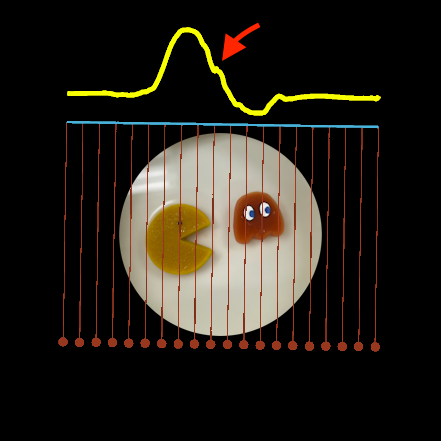}
    \caption*{B}
    \label{fig:profile4b}
  \end{subfigure}
    \hfill
  \begin{subfigure}[b]{0.27\textwidth}
    \includegraphics[width=\textwidth]{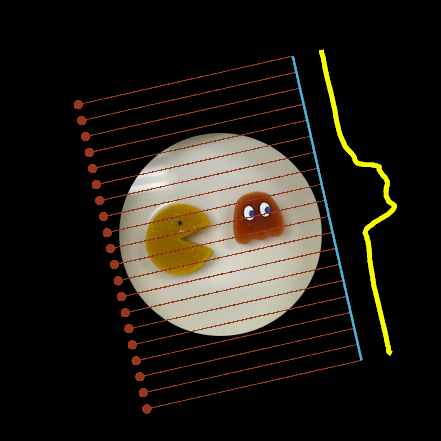}
    \caption*{C}
    \label{fig:profile4c}
  \end{subfigure}
  \caption{Sinogram profiles, superimposed on cropped photos from data collection, with illustrations of virtual X-ray beams. A sampling of  virtual radiographs for targets I, II, III from various directions (A), (B), (C) are plotted in yellow. These are the individual columns of the VHPT sinogram, representing the attenuation of virtual parallel X-ray beams (in red). On III(B) a red  arrow indicates a sharp ``notch'' in the radiograph corresponding to the open mouth of the Pac-Man.}
    \label{fig:profiles}
\end{figure}

In Figure~\ref{fig:profiles} we provide a sampling of virtual radiographs (columns of the sinogram) corresponding to virtual measurement angles (A),(B),(C). Resistive yellow inclusions correspond to peaks in the radiographs, while conductive red inclusions result in valleys.  Observe that the open mouth of the Pac-Man in target III corresponds to a sharp ``notch'' in the radiograph, most visible for angle (B), although this sharp feature is smoothed out in the regularized TV reconstruction. This simple example provides a glimpse of the practical advantages possible with this decomposition of  the EIT problem: by removing some ill-posedness from the inversion process, virtual radiographs can reveal important geometric information that may be lost during the full reconstruction process.

\section*{Discussion}

For five decades now, CT has been a standard technique for uncovering the internal structure of objects, due to its straightforward linear geometry and relatively mild ill-posedness. The EIT  problem is much more challenging than that of CT, and these difficulties lead to characteristic blur in EIT images. However, EIT has advantages in its portability, low cost, and use of nonionizing measurement energy; and EIT has shown promise  for example in the evaluation of lung function in cystic fibrosis patients \cite{mueller2022evaluation}, treatment and monitoring of COVID-19 patients \cite{Jonkman2023,Perier2020,van2020electrical}, classification of stroke \cite{Agnelli2020}, and monitoring fluid flows in chemical and process engineering \cite{powell2008, stephenson2008}. 

Nachman’s 1996 study \cite{nachman1996global} paved the way for the now well-established CGO-based D-bar reconstruction method for EIT \cite{siltanen2000implementation, isaacson2004reconstructions,  mueller2012linear}, wherein the characteristic blurriness of EIT images has been traced back to a nonlinear Fourier-domain viewpoint \cite{knudsen2009regularized}. The exponential behavior of CGO solutions is used to define a  nonlinear Fourier transform for 2D EIT, which is subjected to low-pass filtering at a cutoff frequency dependent on the signal-to-noise ratio, offering regularization but causing blur. Strategies for improving D-bar reconstructions have included the removal of nonlinear blur in a post-processing step \cite{Hamilton2018}, introduction of a Schur-complement based prior \cite{santos2020}, and the inclusion of spatial priors  into the nonlinear Fourier domain \cite{alsaker2016d}. However, in all of these approaches, ill-posedness and nonlinearity were considered to be fundamentally linked.

In this work we invoke related ideas of custom Fourier analysis to solve the EIT problem and yield high-quality reconstructions, but we do so in a way which decouples nonlinearity and ill-posedness. One of the main benefits of the proposed approach, therefore, is enhanced interpretability. The complicated EIT problem is decomposed into a series of comparatively simple modules which isolate the most troublesome features. Ill-posedness is contained entirely within two well-understood modules, permitting application of standard linear techniques. Only the well-posed nonlinear step (D) in Figure~\ref{fig:VHEDprocess2} truly necessitates black-box machine learning, and this module permits clear geometric interpretation. All other nonlinear steps may be treated using explicit analytic methods. Furthermore, we anticipate scholars and engineers will find independent uses for intermediate outcomes within the imaging chain (e.g.\ as inputs to machine learning algorithms), as has already been demonstrated in \cite{Agnelli2020}. 

\subsection*{New Directions} On top of conceptual clarity, our decomposition of the EIT problem unlocks numerous future  imaging possibilities. Via the (deblurred) parallel-beam sinogram, the plethora of reconstruction methods published for linear X-ray tomography in the preceding decades becomes suddenly relevant for nonlinear EIT. Furthermore, while a conductivity image may be readily computed from the sinogram if needed, important geometric information is present in individual columns of the sinogram. These can be interpreted directly if desired, thus removing some  ill-posedness from the imaging process. These ``virtual radiographs'' can for example reveal certain sharp features that may be lost in the regularized conductivity reconstruction, as we have demonstrated in this work.   

In industrial process control for instance, one could generate virtual radiographs to monitor multiphase flow in pipes without X-ray radiation, which in current set-ups necessitates highly specialized facilities, expensive equipment, and radiation protection\cite{powell2008}. The proposed method could greatly reduce the complexity and size of tomographic measuring devices and eliminate risks associated with ionizing radiation. As one example, consider the problem of monitoring the flow of a mixture of opaque fluid and other substances inside a pipeline, e.g.\ in petroleum platforms where long risers transport a mixture of gas, petroleum, and water from the ocean floor. An inexpensive belt of electrodes attached around the pipeline would enable EIT measurements. By approximating the measurement domain as a 2D cross-section of the target, we may compute the VHPT sinogram repeatedly at an appropriate frame rate and plot the time traces of one or more fixed columns as rows in a video frame. Such rolling ``virtual fluoroscopy'' would offer an ``X-ray'' view of inhomogeneities traveling past the 2D measurement domain, without any actual X-rays.

We further anticipate the conceptual and practical advances developed here will extend analogously to other ill-posed problems, including diffuse optical tomography and electrical capacitance tomography. Such modalities have much in common with EIT, including  the use of nonlinear wave phenomena as measurement energy  and the recovery of an unknown PDE coefficient as image reconstruction. These possibilities, so far unexplored, may present new directions for advancement in many fields. CGO solutions have been fruitful in inverse problems related to the Helmholtz and heat equations, and we anticipate VHPT-type methods to be possible for those as well.

\section*{Data Availability and Reproducibility} 
All data and code necessary to reproduce the VHPT reconstructions provided in this work are available in \cite{VHPT_code}. 

\section*{Acknowledgements}
The authors thank S.J.A.~Cowan, N.~Linthacum, T.M.~McKenzie, and M.B.P.~Wong for computational assistance, and also A.V.~Pigatto and N.~Barbosa~da~Rosa~Jr.\ for assistance with data collection. M.A. acknowledges support from the Alphonse A. \& Geraldine F. Arnold Endowment. S.R. and S.S. acknowledge support from the Research Council of Finland (353097). J.P.A. acknowledges support from Agencia Nacional de Promoci\'on Cient\'fica y Tecnol\'ogica (ANPCyT) project PICT 2021-0188. R.M. acknowledges support from the Jane and Aatos Erkko Foundation. 
M.L. acknowledges support from AdG project 101097198 of the European Research Council. Views and opinions expressed are those of the authors only and do not necessarily reflect those of the E.U.

\clearpage
\section*{Supplementary materials}
In this Appendix we elaborate upon some  technical and mathematical details behind the VHPT techniques described in the main text of our work. This material is considered optional reading for those readers who desire deeper understanding of VHPT. 

\subsection*{Details of Computing the DN Matrix} 
\label{sec:realDataToDNmatrix}
Here, we provide details of our calibration approach for completing step (A) of the VHPT imaging chain described in the main text. This involves  computing a DN matrix $\bL_{\sigma}$ 
from   real-world EIT measurements, and our approach  may be thought of as a  practical implementation of the ideas presented in~\cite{garde2021mimicking}.  
To compute a numerical solution of the boundary integral equation
\begin{equation} \label{eqn:BIE}
{M_{\pm \mu}(\,\cdot\,,k)|_{\partial\Omega} + 1 = ({\mathcal P}_{\pm\mu}^k + {\mathcal P}_0)M_{\pm \mu}(\,\cdot\,,k)|_{\partial\Omega},}
\end{equation}
described in the main text, we must first compute a matrix approximation $\bL_{\sigma}$ of the DN map $\Lambda_{\sigma}$ associated with the continuum model
\begin{align}
\nabla\cdot\sigma\nabla u &= 0 \quad \mbox{ in }\Om, \label{eq:cond_eq}\\
 u &= f \quad \mbox{ on }\DOm \label{eq:Dirichlet_cond}.
\end{align}
In practical real-world computations,  using $L$ electrodes arranged around $\partial\Omega$, we measure the voltage distribution arising from each of $L-1$ linearly independent trigonometric current patterns, applied with peak amplitude $I_0$.  Raw current measurements are rearranged into a  matrix $J$ of size $L \times (L-1)$, where the $n$th column corresponds to the basis functions 
\begin{equation} \label{eq:trig-cp}
\phi_n(\theta_\ell)= 
\begin{cases} 
\displaystyle I_0 \cos\left(\frac{(n+1)}{2}\theta_\ell \right) &\quad \mbox{odd } n, \\[3ex]
\displaystyle I_0 \sin\left(\frac{n}{2}  \theta_\ell \right) &\quad \mbox{even } n,
\end{cases}
\end{equation}
$n = 1, \dots, L-1$, where $\theta_\ell, \; \ell = 1, \dots, L$ are taken to be electrode centers.

Corresponding voltage matrices $V_\mathrm{clb}$ and $V_\mathrm{trg}$ are formed for the calibration data (e.g. a homogeneous domain) and the target data, respectively. From these current and voltage measurements we form a discrete matrix approximation $\bL_{\sigma}$ to the DN map $\Lambda_{\sigma}$ via a novel calibration approach which we describe in detail later in this section. 

If instead of considering the Dirichlet boundary condition~\eqref{eq:Dirichlet_cond} to the equation~\eqref{eq:cond_eq} we consider the Neumann boundary condition  
\begin{equation*}
\sigma\frac{\dd u}{\dd \nu} = \phi \quad \text{on } \DOm
\end{equation*}
satisfying $\int_{\DOm} u \, ds = 0$, then the Neumann-to-Dirichlet (ND) map
$\cR_\sigma :H^{-1/2}_{\diamond}(\partial \Omega) \to H^{1/2}_{\diamond} (\DOm)$
is given by
\begin{equation}\label{eq:ND-map}
\cR_{\sigma}: \phi  \longmapsto  u_{|_{\DOm}}.
\end{equation}
Here $H^{s}_{\diamond}(\DOm)$ spaces consist of functions in the standard Sobolev space $H^{s}(\DOm)$ having  mean value zero. 
Given a set of basis functions $\{\phi_n\}_{n>0}$ on $\DOm$, an infinite discrete matrix approximation $\mbR_{\sigma}$ to the ND map~\eqref{eq:ND-map} can be constructed by first computing $\mathcal{R}_\sigma \phi_n$,  $n>0$ and then evaluating the inner products
\[
\langle \phi_m,\mathcal{R}_\sigma \phi_n\rangle
=: ( \mathbf{R}_{\sigma} )_{mn},
\]
where $\langle \cdot,\cdot \rangle $ denotes the usual inner product in $L^{2}(\DOm)$. 
Inverting this matrix $\mathbf{R}_{\sigma}$, one has a matrix approximation $\bL_{\sigma}$ of the DN map $\Lambda_{\sigma}$ restricted to the subspace of zero-mean functions.

However, in practical applications instead of continuum data one has a finite number of voltage and current measurements. 
The most accurate and widely used mathematical model for the electrode measurements is the \emph{complete electrode model} (CEM)~\cite{somersalo1992existence}. Let $e_\ell \subset \DOm$, $\ell=1,\dots,L$ be disjoint open paths modeling the electrodes. The CEM is defined by the conductivity equation~\eqref{eq:cond_eq} and the following boundary conditions: 
\begin{align}
\label{eq:cem1}
z_\ell \,  \sigma \frac{\dd u}{\dd \nu} + u  &= V_\ell, \quad \,  \text{on  } e_{\ell}, \;  \ell=1,\ldots, L   \\
\label{eq:cem2}
\frac{\dd u}{\dd \nu} &= 0, \quad \;  \text{ on }  \DOm \setminus \cup_{\ell = 1}^{L}e_\ell  \\
\label{eq:cem3}
\int_{e_\ell}  \sigma \frac{\dd u}{\dd \nu} \dif s &= I_\ell , \quad \, \,  \ell=1,\dots,L,
\end{align}
where $V_\ell$ is the constant electric potential on electrode $e_\ell$, $z_\ell$ is the contact impedance at electrode $e_\ell$, and $I_\ell$ denotes the net current on electrode $e_\ell$. The \emph{absolute} measurements of the CEM are modeled by the electrode current-to-voltage map
$\cR_\sigma^{\text{cem}}:  \R^L_{\diamond} \to  \R^L_{\diamond}$ given by 
\begin{equation*}\label{eq:ND-map-CEM}
\cR_{\sigma}^{\text{cem}}: I  \longmapsto  V.
\end{equation*}

Analogously to the formation of the infinite discrete $\mathbf{R}_\sigma$ computed from the continuum model, we can calculate a finite matrix approximation $\mbR^{\text{cem}}_{\sigma}$ of the electrode current-to-voltage map $\cR_{\sigma}^{\text{cem}}$ and then invert this matrix to compute a matrix approximation $\bL^{\text{cem}}_{\sigma}$ of the DN map $\Lambda_{\sigma}$. 
However, in this study we follow ideas presented in~\cite{garde2021mimicking} and use CEM data to form an  approximation $\mathbf{Y}^{\text{cem}}_\sigma$ to the \emph{relative} continuum ND map 
\[ 
\Upsilon_\sigma := \cR_\sigma - \cR_1,
\]
where $\cR_1$ is the ND map corresponding to homogeneous conductivity $\sigma \equiv 1$ on $\Omega$. Then, from $\mathbf{Y}^{\text{cem}}_\sigma$ we compute the desired finite DN matrix $\bL_\sigma$.

Since we take  $\Omega$ to be the unit disc, we assume the center of electrode $e_\ell$ is given by $x_\ell = \text{e}^{i\theta_\ell}$ where ${\theta_\ell = \ell \frac{2 \pi}{L}}$, for $\ell=1,\ldots, L$. 
In this study, the continuum current patterns are given by the trigonometric basis functions~\eqref{eq:trig-cp}. 
Under the previous assumptions the numerical algorithm of the method proposed in~\cite{garde2021mimicking} to compute a matrix approximation of the continuum DN map $\Lambda_\sigma$ reduces to the following steps:

\begin{itemize}
	\item \textbf{Step 1:} For even $L$, take $n=1,\ldots,L-2$ and define
    \begin{equation*}\label{eq:I_cp}
	I^{(n)} := \frac{2\pi}{L} 	
	\left[ 
	\phi_n(\theta_{1}), \ldots, \phi_n(\theta_L) 
    \right]^{T}   
	\in  \mathbb{R}^{L}_\diamond.
	\end{equation*}
    We discard $I^{(L-1)}$ since the  CGO-based reconstruction method described in the main text requires an even number of basis functions to ensure invertibility of the  matrix approximation to a tangential derivative operator used in the construction of the projection operators $\mathcal{P}^k_{\pm\mu}$ and $\mathcal{P}_0$ in equation \eqref{eqn:BIE} (see \eqref{eq:tangential_D}, defined later in this Appendix). If $L$ is even (as is the case with most EIT systems), then $L-1$ is odd and we must therefore eliminate this extraneous current pattern. 
			
	\item\textbf{Step 2:} For $n=1,\ldots,L-2$, perform the measurements of CEM with $I^{(n)}$ as input current to retrieve 
    \begin{align*}
       V_{\sigma}^{(n)}&= \cR_{\sigma}^{\text{cem}} I^{(n)} \in \R^{L}_\diamond \\
       V_1^{(n)}&=  \cR_{1}^{\text{cem}}  I^{(n)} \in \R^{L}_\diamond.         
    \end{align*}    
    Then, define $U_{\sigma}^{(n)}= V_{\sigma}^{(n)} - V_1^{(n)}$ 
    and form the matrix $\mathbf{U} = [ U^{(1)} \, U^{(2)} \ldots U^{(L-2)}].$
		
	\item\textbf{Step 3:} Construct the matrix $\mathbf{I} = [ I^{(1)} \, I^{(2)} \ldots I^{(L-2)}].$ 
	
	\item\textbf{Step 4:} Compute a matrix approximation of the \emph{relative} electrode current-to-voltage map as
	\[
	\mathbf{Y}^{\text{cem}}_\sigma =    \mathbf{I}^T \mathbf{U}. 
	\]

    \item\textbf{Step 5:}  Finally, a matrix approximation of the DN map $\Lambda_{\sigma}$ acting on the the basis~\eqref{eq:trig-cp} is given by  
	\begin{equation*}
	\bL^{\text{cem}}_{\sigma} = (\mathbf{Y}^{\text{cem}}_\sigma + \mathbf{R}_{1})^{-1},
	\end{equation*}
where $\mathbf{R}_{1}$ is an ideal reference ND matrix given by
\begin{equation*}\label{eq:ND_ideal}
	\mathbf{R}_{1} = \text{diag}(1,1,\frac12,\frac12,\ldots,\frac{2}{L-2},\frac{2}{L-2}).
 \end{equation*}
	
\end{itemize}
This process  may be thought of as a form of difference imaging, in that to compute the approximate DN map $\bL^{\text{cem}}_{\sigma}$ we use the difference measurements  $\mathbf{U} = \mathbf{V}_{\sigma} -\mathbf{V}_{1}$, as opposed to the absolute imaging scenario wherein we would compute and invert the finite matrix approximation $\mbR^{\text{cem}}_{\sigma}$ of the current-to-voltage map $\cR_{\sigma}^{\text{cem}}$. In practical computations with experimental data, EIT difference imaging has been demonstrated to provide significant advantages over absolute imaging in terms of robustness to modeling errors and measurement noise~\cite{brazey2022robust}.

To compute the voltage matrix $\mathbf{U}$ in Step 2 one needs a matrix approximation of $\cR_1^{\text{cem}}$. Therefore, to directly implement the method proposed in~\cite{garde2021mimicking} using real-world data, one needs EIT measurements corresponding to homogeneous unit conductivity, that is $\sigma=1$, but this is not always possible. However, in laboratory experiments we may collect measurements on an ``empty'' tank filled with homogeneous saline. Here, we present a novel calibration step that allows us to approximate $\cR_1$ from such empty-tank measurements. Then, applying the same calibration step to the measurements collected on the target, we obtain an approximation of $\cR_{\sigma}$, which subsequently yields a good matrix approximation of $\Lambda_{\sigma}$ in the spirit of the method proposed in~\cite{garde2021mimicking}. The procedure consists of the following steps:

\begin{itemize}

    \item\textbf{Step (i):}  Collect EIT measurements on a tank filled with homogeneous saline solution having (assumed unknown) conductivity $\sigma_{\text{clb}}$, using applied  currents $I^{(n)}$ to obtain corresponding calibration voltages $V_{\text{clb}}^{(n)}$. Then,  form the matrix $\mbR_{\text{clb}}=  \mathbf{I}^T \mathbf{V}_{\text{clb}}$,  where
    $\mathbf{V}_{\text{clb}} = [ V_{\text{clb}}^{(1)} \, V_{\text{clb}}^{(2)} \ldots V_{\text{clb}}^{(L-2)}].$

    \item\textbf{Step (ii):}  Compute the calibration matrix $\mathbf{C}$ such that
	\[
	\mathbf{C}\mbR_{\text{clb}} = \mbR_{1}^{\text{cem}}. \]
	This matrix is given by $\mathbf{C} = \mathbf{R}_{1}^{\text{cem}} (\mbR_{\text{clb}})^{-1}$, where
    $\mathbf{R}_{1}^{\text{cem}} = \mathbf{I}^T \mathbf{V}_1$,  
    and
    $\mathbf{V}_1 = [ V_1^{(1)} \, V_1^{(2)} \ldots V_1^{(L-2)}]$ is a matrix of numerically simulated boundary voltage measurements corresponding to $\sigma \equiv 1$  in $\Omega$, computed by forward-solving  \eqref{eq:cond_eq} subject to \eqref{eq:cem1}-\eqref{eq:cem3} using finite element methods (FEM).

   \item\textbf{Step (iii):}  Perform real-world electrical measurements at the boundary of the target object using applied currents $I^{(n)}$ to obtain the corresponding target voltages $V_{\text{trg}}^{(n)}$. Then, form the matrix 
   $\mbR_{\text{trg}}=  \mathbf{I}^T \mathbf{V}_{\text{trg}}$,  where
    $\mathbf{V}_{\text{trg}} = [ V_{\text{trg}}^{(1)} \, V_{\text{trg}}^{(2)} \ldots V_{\text{trg}}^{(L-2)}],$
   and apply the calibration matrix $\mathbf{C}$ to obtain the modified version

      \[  \mathbf{C} \mbR_{\text{trg}}  
    = \mathbf{R}_{1}^{\text{cem}} (\mbR_{\text{clb}})^{-1}  \mbR_{\text{trg}} . 
    \]  
    
    \item\textbf{Step (iv):}  Compute the matrix approximation of the \emph{relative} electrode current-to-voltage map as
    \[ 
    \mathbf{Y}_{\text{trg}} 
    = \mathbf{C}(\mathbf{R}_{\text{trg}} -  \mathbf{R}_{\text{clb}})   
    = \mathbf{C} \mbR_{\text{trg}}  - \mbR_{1}^{\text{cem}} .
    \]

    \item\textbf{Step (v):}  Finally, a (non-symmetric) matrix approximation of the DN map $\Lambda_{\sigma} $  is given by
    \[     
    \tilde\bL_{\sigma}  
    = ( \mathbf{Y}_{\text{trg}} + \mbR_1 )^{-1}
    =
    ( \mathbf{C} \mbR_{\text{trg}} - \mbR_1^{\text{cem}} + \mbR_1 )^{-1}.
    \]    
    However, since $\Lambda_\sigma$ is a symmetric linear operator, we enforce this condition by defining
    \begin{equation*} \label{eq:DN_approx_GH}
    \bL_\sigma = \frac12 (\tilde\bL_\sigma + \tilde\bL_\sigma^T).
    \end{equation*}
\end{itemize}

\subsection*{Details on the Projection Operators for the BIE}

The DN map $\Lambda_\sigma$ is needed to compute the projection operators ${\mathcal P}_{\pm\mu}^k$ and $\mathcal{P}_0$  used in \eqref{eqn:BIE}, which in VHPT imaging must be solved twice to obtain the needed CGO solutions $M_{\pm \mu}$. Here, for completeness, we provide the definitions of these projection operators. See \cite{astala2006boundary} for further details. 

Consider solutions $f_{\pm \mu}(x,k) = e^{ikx}(1 + \omega^{\pm}(x,k))$ of the Beltrami equation 
\begin{equation*}\label{eq:Beltrami1}
\overline \p_x f_{\mu}(x,k) = \mu(x)\,\overline{\partial_x f_{\mu}(x,k)}, \quad x \in \mathbb{R}^2
\end{equation*}
where the Beltrami parameter is related to the conductivity via
$\mu(x) = \frac{1-\sigma(x)}{1+\sigma(x)}.$ Let $u$ and $v$ denote the real and imaginary parts of $f_\mu$, respectively. The $\mu$-Hilbert transform $\mathcal{H}_\mu : H^{1/2}(\partial \Omega) \to H^{1/2}(\partial \Omega)$ is then defined by 
\begin{equation*}\label{eq:muHilbert}
\mathcal{H}_\mu u := v,
\end{equation*}
and for real-valued $g \in H^{1/2}(\partial \Omega)$ the operator $\mathcal{H}_\mu$ is directly related to the DN map via the formula 
\begin{equation}\label{eq:tangential_D}
\partial_T \mathcal{H}_\mu g  = \Lambda_\sigma g,
\end{equation}
where $\partial_T$ denotes the tangential derivative along the boundary. Then we may define the projection operators conjugated with multiplications by complex exponential functions:
\begin{equation}\label{eq:projOp}
\mathcal{P}^k_\mu := e^{-ikx}\mathcal{P}_\mu e^{ikx},
\end{equation}
where the un-conjugated operator is given by
$$
\mathcal{P}_\mu g = \frac{1}{2} (I + i \mathcal{H}_\mu)g + \frac{1}{4\pi} \int_{\partial \Omega} g \; ds.
$$
The operator $\mathcal{P}_0$ is the special case of \eqref{eq:projOp} corresponding to the DN map $\Lambda_1$ for the constant-conductivity case $\sigma =1$ in $\Omega$. Using these definitions of $\mathcal{P}^k_\mu$ and $\mathcal{P}_0$, in \cite{astala2006boundary} it is shown that the boundary integral equation \eqref{eqn:BIE} holds. For details of the numerical formulation of the projection operators and numerical solution of \eqref{eqn:BIE}, see \cite{astala2011direct}.

\subsection*{Neural Network Details} \label{sec:NNdetails}

\begin{figure}[t!]
\centering
\input{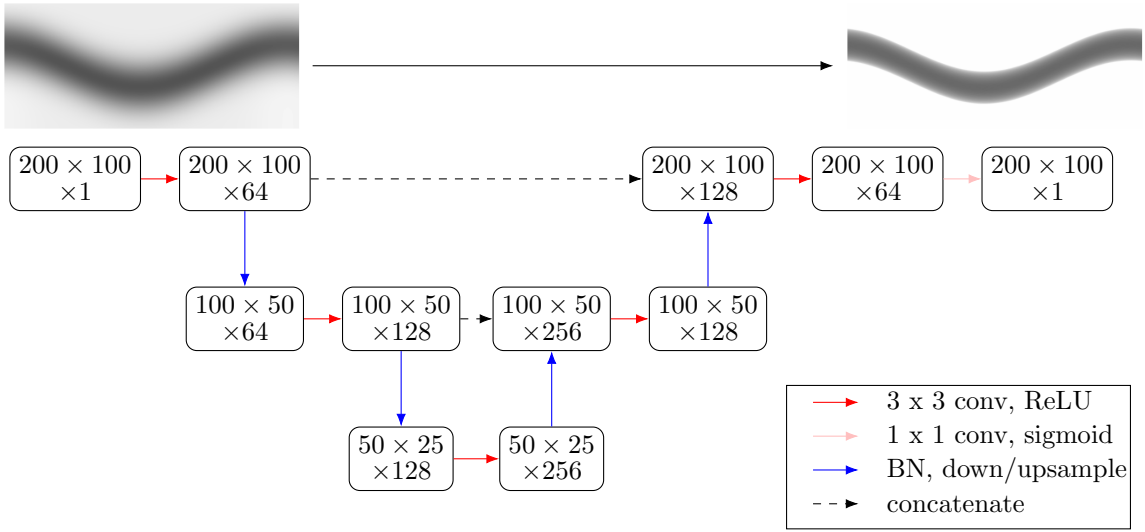}
\caption{Neural network architecture. The size of the feature maps and number of channels is reported for each layer. BN stands for batch normalization.}
\label{fig:nn_arch}
\end{figure}

As described in the main text, we use machine learning to simultaneously complete steps (D) and (E) in the VHPT imaging chain. Our network is a convolutional autoencoder with skip connections, that is, the commonly used U-net architecture \cite{ronneberger2015u}. We used $3 \times 3$ convolution kernels and ReLU activation functions in all layers except the output layer. The final layer used $1 \times 1$ convolution kernels and a sigmoid activation function. Batch normalization was used before performing each downsampling and upsampling, and the encoder and decoder sides of the network were concatenated using skip connections. See Figure \ref{fig:nn_arch} for an illustration of the network architecture.

For training, we used 7000 training pairs of sinograms, from which 1400 pairs were left for validation purposes. We used mean squared error (MSE) as the loss function, and Adam optimization \cite{kingma2014adam} with a learning rate of 0.001. Satisfactory learning results were obtained after 40 epochs of training. 

\subsection*{Details of Reconstruction from a Sharp Sinogram}

Here we elaborate upon the final step (G) of the VHPT process. Given the sharp sinogram $S$ as described in the main text, one may view the columns of $S$ directly to gain geometric information from these ``virtual radiographs,'' thus avoiding some ill-posedness present in the full reconstruction process. If desired however, any standard parallel-beam CT reconstruction algorithm may be applied to transform the sinogram to the spatial imaging domain, thus adding one additional mildly ill-posed step. The simplest option is to use  the traditional filtered back-projection (FBP) formula to reconstruct the Beltrami coefficient $\mu$:
\begin{equation*}
\label{eq:radon_inv}
   \mu(x) = \frac{1}{2}\mathcal{F}^{-1}\Big(b(\xi)\, \mathcal{F}(R^* R\mu)\Big)(x),\quad b(\xi)=|\xi|,
\end{equation*}
where we use $R\mu \approx S$, the back-projection operator (adjoint of the Radon transform) is denoted by $R^*$, $\mathcal{F}$ is the  2D Fourier transform and $\mathcal{F}^{-1}$ its inverse. If a regularized reconstruction is desired, one may use (for example) total variation (TV) regularization \cite{rudin1992nonlinear,mueller2012linear}. Using TV reconstruction, we take $\mu$ to be the solution to the minimization problem
\begin{equation*} \label{min_problem}
 \mu =    \mathop{\mathrm{arg\,min}}_{\mu'} \bigg\{  \frac{1}{2} \| R\mu' - S\|_2^2 + \alpha \textrm{TV}(\mu')\bigg\},
 \end{equation*}
where $R\mu'$ denotes the Radon transform of the parameter $\mu'$,  $\alpha >0$ is the regularization parameter and $\textrm{TV}(\mu')$ is the TV regularization term over the domain $\Omega$, given informally by
\begin{equation*}
    \textrm{TV}(\mu') =  \int_{\Omega} |\nabla \mu'| \dif x.
\end{equation*} 
The regularization parameter values for the reconstructions presented in Figure 2 are the following: $\alpha_{\text{I}} = 0.009$, $\alpha_{\text{II}} = 0.006$, and $\alpha_{\text{III}} = 0.02$. Once we have a reconstruction of the Beltrami coefficient $\mu$ using any desired method, the conductivity $\sigma$ is obtained by the simple nonlinear transformation $\mu :=\frac{1-\sigma}{1+\sigma}.$ This completes the VHPT imaging chain.

\subsection*{Mathematical Background of VHPT}

How is it mathematically possible to  convert the information from measurements of curved electric flow into straight-line tomography? Let us discuss the mathematical framework of the previous work \cite{greenleaf2018propagation}, for those interested in a deeper analysis. 

\subsubsection*{Connecting the Wave Equation and the Conductivity Equation}

Consider the linear wave equation 
$(\frac{1}{c^2}\partial_t^2-\nabla\cdot\nabla)W(x,t)=0$, where $c$ is wave speed and $\nabla=(\partial_{x_1},\partial_{x_2})$. Let $\theta\in\mathbb{R}^2$ be a unit vector. Then the wavefront of a plane wave traveling in direction $\theta$ is given by the solution of the wave equation,
\begin{equation}\label{planewave}
 W(x,t) = \delta(t-\frac {\theta}c\cdot x),
\end{equation}
where $\delta$ is Dirac's delta function. The solution \eqref{planewave} is useful for modeling electromagnetic waves such as X-rays, and below we will connect it to virtual X-rays arising from EIT. 

The main idea of CGO solutions is that their leading-order behavior is exponential, enabling the use of  Fourier transform techniques. For example, the complex exponential $u_0(x)= e^{i\tau \theta\cdot x-\tau\theta^\perp \cdot x}$ solves the constant conductivity equation $\nabla \cdot \nabla u_0=0$. Here $\theta^\perp\in\mathbb{R}^2$ is a unit vector perpendicular to $\theta$. The parameter $\tau\in \mathbb R$ sets both the frequency of spatial   oscillations and the rate of exponential decay (or growth) in a half-plane. Note that we can write $u_0(x)= e^{ik \cdot x}$, where $k=\tau \eta$ is a complexified wave vector with generalized direction $\eta=\theta+i\theta^\perp\in\mathbb{C}^2$. Such non-physical complexification was introduced by Faddeev \cite{Faddeev1966} and later made rigorous in \cite{sylvester1987global}.

When the conductivity is not constant, we write the electric potential as $u(x)=u_0(x)w(x)$, where 
$w=1+w_{sc}$ is the scaled electric potential that satisfies the scaled conductivity 
equation $(\nabla+i\tau \eta)\cdot \sigma (\nabla+i\tau \eta) w=0$. The function $w_{sc}$ can be considered  a perturbation caused by the changes of the conductivity $\sigma(x)$ inside the target. To connect EIT to X-ray imaging, we fix the direction $\theta$ and consider the scaled electric potential $w(x,\tau)$ as a function  of $x\in \mathbb R^2$ and $\tau\in \mathbb R$.  Fourier transformation with respect to the spatial frequency variable $\tau$ yields a function $\widehat w(x,t)$ satisfying the complex principal-type equation
\begin{equation}\label{Cprintype}
(\nabla+\eta\partial_t)\cdot \sigma (\nabla+\eta\partial_t) \widehat w=0. 
\end{equation}
The Fourier variable $t\in \mathbb R$ is called {\it pseudo-time}.

Now in a homogeneous medium, since $\eta\cdot \eta=0$, equation \eqref{Cprintype} takes the form \begin{equation*}\label{Cprintype0}
(\nabla+2\eta \partial_t)\cdot \nabla \widehat w_0(x,t)=0,
\end{equation*}
which allows two interesting solutions. The first is  $\delta(t-2\theta \cdot x)$, analogous to a plane wave of the form \eqref{planewave} traveling in pseudo-spacetime with wave speed $c=1/2$. The second solution $\delta(t-C)$ with constant $C\in\mathbb R$ is parallel to the $x$-plane and it can be considered  a plane wave that travels with an infinite speed. The latter wave is very useful for VHPT because it carries information about conductivity changes to the boundary, which can be detected and processed. In particular as the wave $\delta(t-C)$ travels with infinite speed,
it carries information about the (pseudo-)time at which the incident wave $\delta(t-2\theta \cdot x)$ has
interacted with the discontinuities of the conductivity $\sigma$. Note that the unit vector $\theta$ determines the direction of the virtual X-rays. Any fixed $\theta$ produces one column in the VHPT sinogram. 

\begin{figure}[t!]
\centering
\includegraphics[width = 0.5\textwidth]{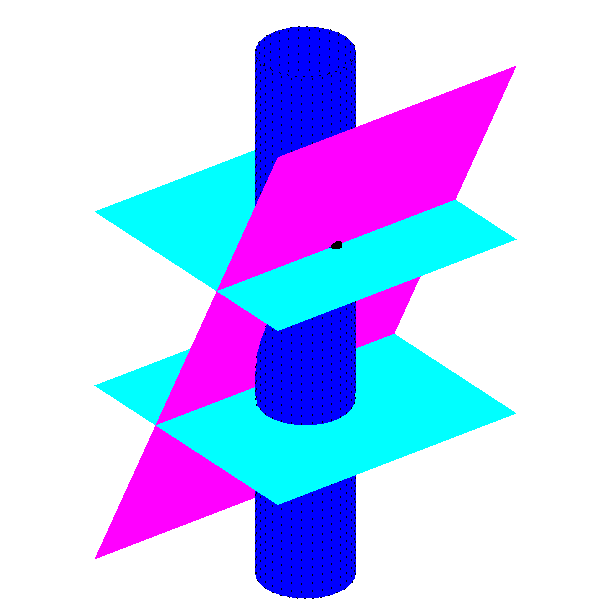}
\caption{
Illustration of the propagation and reflection of singularities in pseudo-spacetime.  Consider the conductivity $\sigma(x)=1+\chi_{B(0,r_0)}(x)$. The dark blue pillar is a vertical cylinder whose base is the singular support of $\mu$, here $\partial B(0,r_0)$. For 
 fixed $\vp$, the map $(x,t)\mapsto \tw_1(x,t,e^{i\vp})$ is singular on  three planes.
The magenta plane is the singular support of the incident wave $\rho_{\vp}(x,t)=c\delta'(t+2\re(e^{i\vp}x))$.
We have $\widehat \w _1(x,t,e^{i\vp})=\dbar^{-1}_x(\mu\,\cdotp \rho_{\vp})$ where $\dbar^{-1}_x$ propagates the singularities of 
the product $\mu\,\cdotp \rho_{\vp}$ in the light-blue horizontal planes $t=2r_0$ and 
$t=-2r_0$.}
\label{fig:singularities}
 \end{figure}

Multiple scattering of virtual waves in pseudo-spacetime produces complicated features in the total wave $\widehat w$. However, rigorous analysis of multiple reflections shows that the reflected waves propagate only along planes of the forms $\delta(t-2\theta \cdot x-C)$ or $\delta(t-C)$. While this results in an infinite ladder of signals scattered multiple times, the structure is still so systematic that we can extract first-order scattering using a relatively shallow neural network. (Curiously,  rectifying the effect of multiple reflections is more complicated for the classical wave equation but can be done using the scattering control method \cite{Caday2019}.) See Figure~\ref{fig:singularities} for an illustration of the reflection and propagation of singularities in pseudo-spacetime. 

\subsubsection*{Connecting the Radon Transform and 
the Conductivity Equation}
\textcolor{magenta}{}
The propagation of certain types 
of oscillations in the solutions of the conductivity 
equation is analyzed in detail in \cite{greenleaf2018propagation}
using microlocal analysis and so-called complex principal-type
operators, introduced by Duistermaat and H\"ormander \cite{duistermaat1972fourier}. More precisely, it is shown in Greenleaf {\it et al.} \cite{greenleaf2018propagation} that the conductivity 
equation may be transformed to an equation of which the principal part is
the complex derivative $\overline \partial_x$,  by conjugating the equation with (i) operators of multiplication by a complex exponential function, and (ii) appropriate Fourier integral operators. This implies that the oscillations in the solutions
of the conductivity equation propagate along two-dimensional spaces in the suitable phase space and satisfy a relation involving the two-dimensional Radon transform. Below we will derive the relationship between the Radon transform and 
 the conductivity equation in a simpler way, using only elementary complex analysis and the properties of the Fourier transform.
To introduce the complex formulation of the conductivity equation
we identify the real-valued coordinates $x=(x_1,x_2)\in \R^2$ with complex-valued points $x=x_1+ix_2\in \C$. We recall that the complex Wirtinger derivatives are
$$ 
\overline \partial _x=\frac 12( \frac{\p}{\p x_1}+i\frac{\p}{\p x_2}),\quad  \partial _x=\frac 12( \frac{\p}{\p x_1}-i\frac{\p}{\p x_2}),\quad x\in \C.
$$
Let $k\in \C$ be the complex wave number, factored as $k=\tau \theta$ where $\tau\in \R$ and $\theta\in \C$, $|\theta|=1$. We call $\tau$ the complex frequency, and later we will introduce 
the pseudo-time variable $t$, which is the Fourier-domain variable
corresponding to $\tau$.

The solution $u$ of the conductivity equation \eqref{eq:cond_eq}
can be written
as
\beq\label{u representation}
u(x)=\Re\,(f_{\mu}(x)) 
\eeq
where $f_{\mu}(x)$ 
is the solution of the Beltrami equation 
\begin{equation*}\label{eq:Beltrami2}
\overline \p_x f_{\mu}(x) = \mu(x)\,\overline{\partial_x f_{\mu}(x)}, \quad x \in \mathbb{C}.
\end{equation*}
Until otherwise stated, we treat $x$ as a complex number $x=x_1+ix_2$.
The representation \eqref{u representation} is analogous to the fact that
any harmonic function can be written as
the real part of an analytic function.

Moreover, Astala and P\"aiv\"arinta \cite{astala2006boundary} constructed unique complex geometrical optics (CGO) 
solutions  
for  
the Beltrami equation
with specific asymptotics at infinity,
\begin{align}
\overline \p_x f_\mu(x,k) &= \mu(x)\,\overline{\partial_x f_\mu(x,k)},\quad 
  x\in \C,  \label{eq:Beltrami3} \\
  f_\mu(x,k) &= e^{ikx}(1+{\mathcal O}(|x|^{-1})) \hbox{ as }|x|\to \infty,  
\end{align}
which depend on 
${x}$ and a complex wave number $k\in \C$.
When $\sigma({x})$ is equal to 1 in the complement of the domain $\Omega$
and  $0<c_1\le \sigma({x})\le c_2$,
the Beltrami coefficient
$\mu$ is   supported in $\Omega$ and  satisfies $\|\mu\|_{L^\infty(\C)}<1$. 
Using such solutions of the Beltrami equations, we will below consider the function 
\beq\label{u representation2}
u(x,k)=\Re\,(f_\mu(x,k))
\eeq
which satisfies $\nabla_{x} \,\cdotp \sigma({x}) \nabla_{x} u({x},k)=0$. It is shown in \cite{astala2006boundary} that
the map $\Lambda_\sigma:u|_{\p \Omega}\to  n\cdot \sigma \nabla u|_{\p \Omega}$ determines the function $f_{\mu}(x,k)$
for $x\in \C\setminus \Omega$ and $k\in \C$.
Let us write the CGO solutions of  
\eqref{eq:Beltrami3} in the form 
\ba
& &{f_\mu(x,k)=e^{ikx}\w (x,k),}\\
& &\w (x,k)=1+{ \w_{sc} (x,k),}\\
& &\w_{sc} (x,k)=
{\mathcal O}(|x|^{-1}),\quad\hbox{as }x\to \infty.
\ea
 Then, \eqref{eq:Beltrami3} yields
 \ba
 \overline \p_x \w _{sc} - e_{-k}\mu  (\overline{\p_x\w _{sc}  - i k \w _{sc}})=-i\overline{k}e_{-k}\mu.
\ea

The solid Cauchy transform $\overline \p_x^{-1}$  is the operator
\beq\label{C-transform}
  \overline \p_x ^{-1} f(x) &=&\frac{1}{\pi}\int_\C \frac{1}{x-x'}f(x')\,d^2 x', 
\eeq
where $d^2x'=dm(x')$ is the Lebesgue measure of $\C=\R^2$.
 Then, we see that $v_{sc}(x,k)$ satisfies
  $$ \big( I-\A\big) \w _{sc}=
  -i\overline k\,  \overline \p_x^{-1} (e_{-k}\mu),$$
where $\A$ is the real-linear operator
$$
\A :v\to \overline \p_x^{-1} (  e_{k}\mu  \rho (  \p_x+ik   )  v)
$$
and $\rho(f):=\overline{f}$ denotes complex conjugation.
Using the equation $ \big( I- \A\big) \w _{sc}=
 -i\overline k\,  \overline \p_x^{-1} (e_{-k}\mu)$ and 
 the Neumann series for $( I- \A)^{-1}$,
we can write formally (the discussion on the precise sense in which this series converges
is omitted in this paper, see \cite{greenleaf2018propagation}
 for detailed analysis)
\beq\label{scattering series A}
 \w _{sc}\sim \sum_{n=1}^\infty \w _n,\quad \w _1= -i\overline k \, \overline \p_x^{-1} (e_{-k}\mu),\quad
  \w _{n+1}=\A \w _n.
 \eeq
Next we consider the term corresponding to {single scattering,}
\ba
{\w _1=  -i\overline k \, \overline \p_x^{-1} (e_{-k}\mu).}\nonumber 
\ea
The Dirichlet-to-Neumann map $\Lambda_\sigma$  determines 
$\w_{sc} (x,k)$  for ${x\in\partial\Omega}$. However, below we consider
the properties of the first term $\w _1(x,k)|_{x\in\partial\Omega}$
in the series \eqref{scattering series A}.
The higher-order terms $\w _n(x,k)|_{x\in\partial\Omega,}\, n\ge2$ contribute ``multiple scattering,'' which
explain artifacts in numerics.

Next, we  analyze the term $\w _1(x,k)$ corresponding to the single-order scattering in the pseudo-time domain.
To do this we decompose the complex wave number $k$ as  $k=\tau \theta,$ 
where $\theta=e^{i\varphi }$. We use the Fourier transform of the function $w (x,\tau)$,
 that is,
$$
\widehat w(x,t)=\F_{\tau\to t} w (x,t)=\int_\R e^{-i t\tau }w (x,\tau)\,d\tau.
$$
For a function $f(x,t)$ that depends on the $x$ and $t$ variables, the formula \eqref{C-transform} implies
\begin{align*}
  \overline \p_x ^{-1} f(x,t)
  &= \frac{1}{\pi} \int_\C \frac{1}{x-x'}f(x',t)\,d^2 x'\\
  &= \frac{1}{\pi}\int_\R \int_\C \delta(t-t')\frac{1}{x-x'}f(x',t')\,d^2 x'dt'.
\end{align*}
This means that $ f(x,t)\to \overline \p_x ^{-1} f(x,t)$
operates in each  $t=$constant plane independently of
the values of $f(x,t)$ in the other planes.

Using the polar coordinates  $k=\tau e^{i\vp}$,
we see for $w=w(x,\tau),$
$$
e_{-k}w=e^{-2i\re(kx) }w(x,\tau)=
e^{-2i\tau \re( e^{i\vp}x) }w(x,\tau).
$$
Thus, it holds for the Fourier transform that
\ba
\widehat {\e_{-k}w} (x,t) &=&\F_{\tau\to t} (e_{-k}w) (x,t)\\
&=&\int_\R e^{-i t\tau }
e^{-2i\tau \re( e^{i\vp}x) }w (x,\tau)\,d\tau
\\
&=&\widehat w(x,  t+2\re( e^{i\vp}x)).
\ea
This can we viewed as `tilting' the
$\{t=0\}$ plane in the space  $\C\times \R$, having the coordinates $(x,\tau)$, by using 
the operation $(x,\tau)\to (x,t+2\re( e^{i\vp}x)).$

As the Fourier transform of
$\mu(x,\tau)=\mu(x)=\mu(x){\bf 1}(\tau)$ is 
$\widehat \mu(x,t)=\mu(x)\delta(\tau)$, we see that
\ba
\F_{\tau\to t} (e_{-k}\mu) (x,t)
&=&\delta( t+2\re( e^{i\vp}x))\mu(x).
\ea
Thus,
the Fourier transform  of the single scattering term
$$
\w _1(x,t,e^{i\vp})= -i\overline k \, \overline \p_x^{-1} (e_{-k}\mu)$$
is equal to
\ba
\widehat \w _1(x,t,e^{i\vp})
=2e^{-i\vp}\int_{\C}\frac 1{x-x'}\, {\delta'(t+2\re(e^{i\vp}x'))}\mu(x')\, d^2x'.\nonumber
\ea
Due to this, we define the linear operator $T_1$ that maps
a compactly supported function $\mu(x')$ in $\Omega$ to the corresponding single scattering term $\widehat v_1$,
$$
\mu(x')\mapsto (T_1\mu)(x,t,e^{i\vp})=\tw_1(x,t,e^{i\vp}).
$$
The Schwartz kernel of $T_1$ is
\ba
K_1(x,t,e^{i\vp},x')=\Big(\frac{2e^{-i\vp}}{x-x'}\Big)\delta'(t+2\re(e^{i\vp}x')),\nonumber
\ea
that is,
\ba
(T_1\mu)(x,t,e^{i\vp})=\frac {\p}{\p t}\int_\C \frac{2e^{-i\vp}}{x-x'}\delta(t+2\re(e^{i\vp}x'))\mu(x')d^2x'.
\ea

Next, we aim to introduce a filtered back-projection formula.
To this end, 
let us consider the `complex average' of $\tw_1(x,t,e^{i\vp})$ over $x\in \p \Omega$,
\ba
\tw_1^a(t,e^{i\vp}) = \frac1{2\pi i}\int_{\partial\Omega} \tw_1(x,t,e^{i\vp}) { dx.}\nonumber
\ea
 Let $T_1^a$ be the linear operator $T_1^a:\mu\mapsto\tw_1^a$.  That is,
 \ba
(T_1^a \mu)(t,e^{i\vp})=\frac {\p}{\p t}\int_\C  2e^{-i\vp} \delta(t+2\re(e^{i\vp}x'))\mu(x')d^2x',
\ea
which we obtain by using the Cauchy theorem and the fact that that $x'\in \Omega$, 
so that 
$$
\frac1{2\pi i}\int_{\partial\Omega} \frac{2}{x-x'}dx
=\frac1{2\pi i}\int_{\partial\Omega} \frac{2}{\frac{x}{x'}-1}\frac{1}{x'}dx
=\frac1{2\pi i}\int_{\{x''=x/x':\ x\in \partial\Omega\}} \frac{2}{x''-1}dx''=2.
$$
We observe that $T_1^a$ is a generalized Radon transform:
When we denote 
$n_{\vp}=(\cos({\vp}),\sin({\vp}))\in \R^2$  we have
\beq
\re(e^{i\vp}x)=n_{-\vp}\cdot {x},
\eeq
with $x=x_1+ix_2$ in the left hand side and $x=(x_1,x_2)$ in the right hand side. In what follows, we continue to consider $x$ as a vector in $\R^2$. Then
$$
(T_1^a \mu)(t,e^{i\vp})
=2e^{-i\vp} \frac {\p}{\p t}\int_{\R^2} \delta(t+2n_{-\vp}\cdot { x})\mu({ x})d{ x}
=e^{-i\vp} \frac {\p}{\p t}{\mathcal R}\mu({-\vp},\frac 12 t),
$$
where 
$$
{\mathcal R}\mu({\vp},s)=\int_{n_{\vp}\cdot { x}=s}\mu({ x})d\ell({ x})
$$
is the 2-dimensional Radon transform of the function $\mu$
and $d\ell$ is the length measure on a line.

The filtered back-projection formula in \cite{natterer1986computerized} for the inverse Radon transform leads to the following theorem (for a detailed proof, see \cite{greenleaf2018propagation}).
\begin{theorem} The operator
 $(T_1^a)^* T_1^a$ is a pseudo-differential operator of order $1$ with
$$\sigma_{\mbox{\tiny prin}}((T_1^a)^* T_1^a)(x,\zeta)= 2\pi^2|\zeta|,\quad x\in\Omega,\ \zeta\in \R^2\setminus 0,$$
where $\sigma_{\mbox{\tiny prin}}$ stands for the principal symbol, and 
\beq\label{FBP}
\frac 1{2\pi^2}(-\Delta)^{-\frac12} (T_1^a)^* T_1^a= I+K,
\eeq
where $K$ is a smoothing integral operator 
(a pseudo-differential operator of order $-1$).
\end{theorem}
Formula \eqref{FBP} provides a reconstruction formula for the function $\mu$, and therefore for $\sigma$, which is analogous to the filtered back-projection formula for the Radon transform.


\bibliographystyle{plain}
\bibliography{main}

\end{document}